# ON ASYMPTOTIC STANDARD NORMALITY
# OF THE TWO SAMPLE PIVOT

BY RAJESHWARI MAJUMDAR* AND SUMAN MAJUMDAR

*University of Connecticut and University of Connecticut*

The asymptotic solution to the problem of comparing the means of two heteroscedastic populations, based on two random samples from the populations, hinges on the pivot underpinning the construction of the confidence interval and the test statistic being asymptotically standard Normal, which is known to happen if the two samples are independent and the ratio of the sample sizes converges to a finite positive number. This restriction on the asymptotic behavior of the ratio of the sample sizes carries the risk of rendering the asymptotic justification of the finite sample approximation invalid. It turns out that neither the restriction on the asymptotic behavior of the ratio of the sample sizes nor the assumption of cross sample independence is necessary for the pivotal convergence in question to take place. If the joint distribution of the standardized sample means converges to a spherically symmetric distribution, then that distribution must be bivariate standard Normal (which can happen without the assumption of cross sample independence), and the aforesaid pivotal convergence holds.

**1. Introduction.** The objective of this note is to critically examine the asymptotic solution to the problem of comparing the means of two populations with finite, unequal variances. Let $X_{1,1}, \cdots, X_{1,n_1}$ be a random sample from the first population with mean $\mu_1$ and variance $\sigma_1^2$, and $X_{2,1}, \cdots, X_{2,n_2}$ a random sample from the second population with mean $\mu_2$ and variance $\sigma_2^2$. When all the parameters $\mu_1$, $\mu_2$, $\sigma_1$, and $\sigma_2$ are unknown, under the assumption of *independence of the two samples*, the traditional asymptotic $100(1-\alpha)\%$ confidence interval for $\mu_1 - \mu_2$ is

$$\left(\bar{X}_1 - \bar{X}_2 - z_{\frac{\alpha}{2}}\sqrt{\frac{S_1^2}{n_1} + \frac{S_2^2}{n_2}}, \bar{X}_1 - \bar{X}_2 + z_{\frac{\alpha}{2}}\sqrt{\frac{S_1^2}{n_1} + \frac{S_2^2}{n_2}}\right), \tag{1}$$

where, for $i = 1, 2$, $\bar{X}_i = n_i^{-1}\sum_{j=1}^{n_i} X_{i,j}$ and $S_i^2 = (n_i - 1)^{-1}\sum_{j=1}^{n_i}(X_{i,j} - \bar{X}_i)^2$, and $z_\alpha$ denotes the $100(1-\alpha)^{\text{th}}$ percentile of the standard Normal distribution. The test statistic for testing the null hypothesis $H_0 : \mu_1 - \mu_2 = D_0$ is

$$\frac{\bar{X}_1 - \bar{X}_2 - D_0}{\sqrt{\frac{S_1^2}{n_1} + \frac{S_2^2}{n_2}}}, \tag{2}$$

which is calibrated on the standard Normal scale for the calculation of the observed level of significance and the rejection regions at various levels of significance.

The simple idea underlying this widely used method is frequently presented in undergraduate (see, for example, Wackerly, Mendenhall, and Scheaffer (2008)) and beginning graduate (see, for example, Casella and Berger (2002)) textbooks on mathematical statistics. As Casella and Berger (2002, page 492) explain, the idea is to obtain a point estimator $\widehat{\theta}_n$ of the parameter of interest $\theta$ with variance $\sigma_n^2$ such that $\left(\widehat{\theta}_n - \theta\right)/\sigma_n$ is asymptotically standard Normal and, if the calculation of $\sigma_n$ involves an unknown parameter other than $\theta$, a consistent estimator $\widehat{\sigma}_n$ of $\sigma_n$, so that $\left(\widehat{\theta} - \theta\right)/\widehat{\sigma}_n$ is also asymptotically standard Normal by Slutsky's theorem [Mukhopadhyay (2000, Theorem 5.3.3)]; $\left(\widehat{\theta} - \theta\right)/\widehat{\sigma}_n$ can then be used as a pivot for $\theta$, and the tools for inference on the Normal mean can be applied. However, Mukhopadhyay (2000) seems to be the only one to explicitly consider the two sample problem from this point of view. Assuming that the two samples are independent and the ratio $n_1/n_2 \to \delta \in (0, \infty)$, he asserts (page 544) the asymptotic standard Normality of the pivot

$$W_{n_1,n_2} = \frac{\bar{X}_1 - \bar{X}_2 - (\mu_1 - \mu_2)}{\sqrt{\frac{S_1^2}{n_1} + \frac{S_2^2}{n_2}}} = \frac{(\bar{X}_1 - \mu_1) - (\bar{X}_2 - \mu_2)}{\sqrt{\frac{S_1^2}{n_1} + \frac{S_2^2}{n_2}}} \tag{3}$$

underpinning the formulas of the confidence interval in (1) and the test statistic in (2).

While asymptotically studying a statistical problem using two samples, many authors (see, among others, DasGupta (2008), Pyke and Shorack (1968), Ramdas, Trillos, and Cuturi (2017), and van der Vaart and Wellner (1996)) require that $n_1/n_2 \to \delta \in (0, \infty)$ as $n_1, n_2 \to \infty$, though this requirement risks rendering the asymptotic justification of the finite sample approximation invalid. For example, the two sample problem with $n_2 > n_1(\ln(\ln n_1))$ (think $n_1 = 50$ and $n_2 \geq 69$) is a fairly common design, at which rate $n_1/n_2 \to 0$; the finite sample approximation involved in using the pivot $W_{n_1,n_2}$ to draw inference on $\mu_1 - \mu_2$ in such a design may not have a rigorous asymptotic justification.

To investigate if we can derive the asymptotic standard Normality of $W_{n_1,n_2}$ without the restriction that the ratio of the sample sizes converges to a finite positive number, we formulate the question in terms of weak convergence. For a separable metric space $\mathcal{S}$, let $\mathcal{B}(\mathcal{S})$ denote the Borel $\sigma$-algebra of $\mathcal{S}$ and $\mathcal{M}(\mathcal{S})$ the set of probability measures on $\mathcal{B}(\mathcal{S})$. Endowed with the topology of weak convergence, $\mathcal{M}(\mathcal{S})$ is metrizable as a separable metric space [Parthasarathy (1967, Theorem II.6.2)]. Since all the random elements under consideration are Borel measurable, convergence in distribution is equivalent to weak convergence of the induced probability measures [van der Vaart and Wellner (1996, page 18)]. With

$$F_{n_1,n_2} \in \mathcal{M}(\Re) \text{ denoting the measure induced by } W_{n_1,n_2} \text{ on } \mathcal{B}(\Re), \tag{4}$$

the asymptotic standard Normality of $W_{n_1,n_2}$ is equivalent to the *convergence* of the double sequence $F_{n_1,n_2}$ to $\Phi$ in $\mathcal{M}(\Re)$, where $\Phi$ denotes the standard Normal measure on $\mathcal{B}(\Re)$. Faced with the question of convergence of a double sequence, we examine both the *iterated* and *double* limits of $F_{n_1,n_2}$.

Before proceeding further, let us spell out without any ambiguity what the assumption of independence of the two samples means in this context. Hereinafter, iid will abbreviate independent and identically distributed. We are making the following three assumptions:

$$\{X_{1,j} : j \geq 1\} - \text{iid sequence of random variables, mean } \mu_1 \text{ and variance } \sigma_1^2 \tag{5}$$

$$\{X_{2,j} : j \geq 1\} - \text{iid sequence of random variables, mean } \mu_2 \text{ and variance } \sigma_2^2 \tag{6}$$

$$\text{the sequence } \{X_{1,j} : j \geq 1\} \text{ is independent of the sequence } \{X_{2,j} : j \geq 1\}. \tag{7}$$

The assumption in (7) is often stated as the $X_{1,j}$'s being independent of the $X_{2,j}$'s. Note that the triplet of assumptions (5), (6), and (7) is equivalent to the pair of assumptions

$$\{(X_{1,j}, X_{2,j}) : j \geq 1\} - \text{iid sequence of random vectors} \tag{8}$$

$$\text{for every } j \geq 1, X_{1,j} \text{ and } X_{2,j} \text{ are independent.} \tag{9}$$

We are not aware of any result on the iterated limits of $F_{n_1,n_2}$. Proposition 1 observes that both iterated limits of $F_{n_1,n_2}$ equal $\Phi$ under only (5) and (6), that is, without (7). As far as the double limit of $F_{n_1,n_2}$ is concerned, the only published result we are aware of is Mukhopadhyay (2000) cited above. However, by requiring that $n_1/n_2 \to \delta \in (0, \infty)$ as $n_1, n_2 \to \infty$, Mukhopadhyay (2000) does not obtain the double limit of $F_{n_1,n_2}$.

We obtain the double limit of $F_{n_1,n_2}$ by using the fact that a double sequence is a net. Let $\mathfrak{N}$ denote the set of natural numbers; then $\mathfrak{N} \times \mathfrak{N}$ is a directed set under the partial ordering $\succeq$ defined by the condition that $(n_1, n_2) \succeq (m_1, m_2)$ if and only if $(n_1 \geq m_1 \text{ and } n_2 \geq m_2)$. A double sequence $\{x_{n_1,n_2} : n_1, n_2 \geq 1\}$ taking values in $\mathcal{S}$ converges to $x \in S$ as $n_1, n_2 \to \infty$ if and only if the corresponding net $\{x_\alpha : \alpha \in \mathfrak{N} \times \mathfrak{N}\}$ converges to $x$. Thus, our objective reduces to obtaining

$$\lim_\alpha F_\alpha = \Phi. \tag{10}$$

Proposition 2 shows that (5) and (6) are not sufficient for (10). Proposition 3 shows that (10) is implied by the convergence of the joint distribution of the standardized sample means to a spherically symmetric distribution, which implies that the limiting spherically symmetric distribution is bivariate standard Normal (see Remark 1). Corollary 1 and Remarks 3 and 4 investigate the question of necessity of the convergence of the joint distribution of the standardized sample means to the bivariate standard Normal distribution for (10), to which Proposition 4 furnishes a partial answer, and Remark 6 outlines the setup wherein a complete answer is obtainable. Detailed statements of these results and their proofs constitute Section 2. The technical results that we draw upon for the proofs of our results are assembled in three Appendices, A.1, A.2, and A.3.

**2. Results and Proofs.** In what follows, unless otherwise specified, we will assume that the index $i$ runs from 1 to 2 and $\{X_{i,j} : j \geq 1\}$ are as in (5) and (6).

**Proposition 1** With $F_{n_1,n_2}$ as in (4),

$$\lim_{n_1 \to \infty} \lim_{n_2 \to \infty} F_{n_1,n_2} = \lim_{n_2 \to \infty} \lim_{n_1 \to \infty} F_{n_1,n_2} = \Phi. \tag{11}$$

For subsequent use, let us introduce the following notations. For $k \in \mathfrak{N}$, let $\overline{X}_{i,k}$ denote $\sum_{j=1}^{k} X_{i,j}/k$ and $S_{i,k}^2$ denote $\sum_{j=1}^{k}(X_{i,j} - \overline{X}_{i,k})^2/(k-1)$. For $\rho \in [-1, 1]$, let $(\Phi \times \Phi)_\rho \in \mathcal{M}(\Re^2)$ denote the bivariate Normal distribution with means 0, variances 1, and correlation coefficient $\rho$. Thus, $(\Phi \times \Phi)_0$ is the product measure $\Phi \times \Phi$, the bivariate standard Normal distribution. For $\theta \geq 0$, let $\mathcal{N}_\theta \in \mathcal{M}(\Re)$ denote the centered Normal measure with variance $\theta$ so that $\mathcal{N}_1 = \Phi$. Let $\xi : \mathfrak{N} \mapsto \mathfrak{N} \times \mathfrak{N}$ denote the order preserving and cofinal map given by $\xi(k) = (k,k)$. Further, let $n_i : \mathfrak{N} \times \mathfrak{N} \mapsto \mathfrak{N}$ denote the order preserving and cofinal map that maps $\alpha \in \mathfrak{N} \times \mathfrak{N}$ to its $i^{\text{th}}$ coordinate $n_i(\alpha)$. For $\alpha \in \mathfrak{N} \times \mathfrak{N}$, let

$$n_\alpha = \min(n_1(\alpha), n_2(\alpha)), \quad m_\alpha = \sqrt{n_1(\alpha)n_2(\alpha)}, \quad \text{and} \quad e_\alpha = n_1(\alpha)/n_2(\alpha).$$

**Proposition 2** Let $\{(X_{1,j}, X_{2,j}) : j \geq 1\}$ be iid, with $(\Phi \times \Phi)_\rho$ being the common distribution. Assume $\rho \neq 0$. Then (5) and (6) hold, but (10) does not.

**Definition** For $k \in \mathfrak{N}$, let $\mathcal{I}(\Re^k)$ denote the set of isometries on $\Re^k$ [Axler (2015, Definition 7.37)]. An element $Q \in \mathcal{M}(\Re^k)$ is called spherically symmetric if the Pettis integral of the identity function with respect to $Q$ is $0 \in \Re^k$, and, for every $\mathcal{O} \in \mathcal{I}(\Re^k)$, $Q \circ \mathcal{O}^{-1} = Q$, where $Q \circ \mathcal{O}^{-1} \in \mathcal{M}(\Re^k)$ is the measure induced by $\mathcal{O}$, that is, $Q \circ \mathcal{O}^{-1}(B) = Q(\{(x_1, \cdots, x_k) \in \Re^k : \mathcal{O}(x_1, \cdots, x_k) \in B\})$ for any $B \in \mathcal{B}(\Re^k)$. Let $\mathcal{S}(\Re^k)$ denote the set of spherically symmetric elements in $\mathcal{M}(\Re^k)$.

**Proposition 3** For $\alpha \in \mathfrak{N} \times \mathfrak{N}$, let $Y_\alpha^{(i)}$ denote the standardized sample mean

$$\sqrt{n_i(\alpha)}(\overline{X}_{i,n_i(\alpha)} - \mu_i)/\sigma_i$$

and $H_\alpha \in \mathcal{M}(\Re^2)$ the measure induced by $Y_\alpha = \left(Y_\alpha^{(1)}, Y_\alpha^{(2)}\right)$. Then (10) is implied by

$$\lim_\alpha H_\alpha = Q \in \mathcal{S}(\Re^2). \tag{12}$$

**Remark 1** The sufficient condition for (10) stated in (12) is equivalent to

$$\lim_\alpha H_\alpha = \Phi \times \Phi. \tag{13}$$

Since $\Phi \times \Phi \in \mathcal{S}(\Re^2)$, (13) implies (12). To show the converse, it suffices to assume (12) and show $Q = \Phi \times \Phi$, equivalently, $Q \circ \langle \,\cdot\, , (a_1, a_2) \rangle^{-1} = \Phi$ for every unit vector $(a_1, a_2)$. Since given any two unit vectors there exists an isometry mapping one to the other, by the spherical symmetry of $Q$, $Q \circ \langle \,\cdot\, , (a_1, a_2) \rangle^{-1} = Q \circ \langle \,\cdot\, , (1, 0) \rangle^{-1}$. Since, by the continuous mapping theorem for nets [Lemma A.7] and (12), $Q \circ \langle \,\cdot\, , (1, 0) \rangle^{-1}$ equals $\lim_\alpha H_\alpha \circ \langle \,\cdot\, , (1, 0) \rangle^{-1}$, it suffices to show that $\lim_\alpha H_\alpha \circ \langle \,\cdot\, , (1, 0) \rangle^{-1} = \Phi$.

For $k \in \mathfrak{N}$, let $P_k^{(i)} \in \mathcal{M}(\Re)$ denote the measure induced by $\sqrt{k}(\overline{X}_{i,k} - \mu_i)/\sigma_i$, which, by (5) or (6), converges to $\Phi$ by the Central Limit Theorem (CLT, hereinafter) [Dudley (1989, Theorem 9.5.6)]. Since convergence in $\mathcal{M}(\Re)$ is metrizable, $n_i$ is order preserving and cofinal, and any subnet of a convergent net has the same limit [Lemma A.2],

$$\lim_\alpha P^{(i)}_{n_i(\alpha)} = \lim_{k \in \mathfrak{N}} P^{(i)}_k = \Phi. \tag{14}$$

Since $H_\alpha \circ \langle \,\cdot\, , (1, 0) \rangle^{-1} = P^{(1)}_{n_1(\alpha)}$, the converse follows from (14). //

**Remark 2** Note that (7) implies

$$H_\alpha = P^{(1)}_{n_1(\alpha)} \times P^{(2)}_{n_2(\alpha)}; \tag{15}$$

since the product of two weakly convergent nets of probability measures converges to the product of the limits [Lemma A.8], (15) implies (13) via (14). Thus, the folklore sufficient conditions (5), (6), and (7), equivalently, (8) and (9), do imply (13), and consequently (10), without requiring $n_1/n_2$ to converge to $\delta \in (0, \infty)$. //

**Remark 3** Does (10) imply (13)? While we cannot construct a counterexample where (10) holds, above and beyond (5) and (6), but (13) does not, we do not believe that the affirmative answer holds because of Corollary 1. //

**Corollary 1** Recall the definition of $W_\alpha$ from (3) and define

$$U_\alpha = \frac{(\overline{X}_{1,n_1(\alpha)} - \mu_1) + (\overline{X}_{2,n_2(\alpha)} - \mu_2)}{\sqrt{\frac{S^2_{1,n_1(\alpha)}}{n_1(\alpha)} + \frac{S^2_{2,n_2(\alpha)}}{n_2(\alpha)}}};$$

let $G_\alpha \in \mathcal{M}(\Re)$ denote the measure induced by $U_\alpha$ on $\mathcal{B}(\Re)$. Then (13) implies

$$\lim_\alpha G_\alpha = \Phi. \tag{16}$$

**Remark 4** If (10) were to imply (13), by Corollary 1, (10) would have to imply (16). We see no reason why (10), under only (5) and (6), would imply (16) (or vice-versa). However, Proposition 4 does connect (10), (13), and (16). //

**Proposition 4** If

$$L = \lim_\alpha H_\alpha \tag{17}$$

exists, and (10) and (16) hold, then (13) holds.

**Remark 5** By (14) and Proposition 9.3.4 of Dudley (1989), $\left\{P_k^{(i)} : k \in \mathfrak{N}\right\}$ is uniformly tight. Consequently, $\left\{P_{n_i(\alpha)}^{(i)} : \alpha \in \mathfrak{N} \times \mathfrak{N}\right\}$, being contained in $\left\{P_k^{(i)} : k \in \mathfrak{N}\right\}$, is uniformly tight as well. By Tychonoff's theorem and Bonneferroni's inequality,

$$\{H_\alpha : \alpha \in \mathfrak{N} \times \mathfrak{N}\} \text{ is uniformly tight.} \tag{18}$$

Can we establish Proposition 4 without assuming (17), substituting it by the conclusion drawn in (18)? We do not think so, though a counterexample eludes us. //

Lemma 1 is a vital cog in the wheel of our investigation of whether (10) implies (13).

**Lemma 1** With $\widehat{F}_\alpha \in \mathcal{M}(\mathfrak{R})$ denoting the measure induced by

$$\widehat{W}_\alpha = \frac{\left(\overline{X}_{1,n_1(\alpha)} - \mu_1\right) - \left(\overline{X}_{2,n_2(\alpha)} - \mu_2\right)}{\sqrt{\frac{\sigma_1^2}{n_1(\alpha)} + \frac{\sigma_2^2}{n_2(\alpha)}}}$$

on $\mathcal{B}(\mathfrak{R})$, (10) is equivalent to

$$\lim_\alpha \widehat{F}_\alpha = \Phi. \tag{19}$$

**Remark 6** Proposition 4 strengthens our belief that (10) does not imply (13), but we are, as mentioned above, unable to construct a counterexample. One of the major obstacles to constructing such a counterexample is the fact that it is simply impossible to get a handle on the asymptotic distribution of either $\{W_\alpha : \alpha \in \mathfrak{N} \times \mathfrak{N}\}$ or $\{Y_\alpha : \alpha \in \mathfrak{N} \times \mathfrak{N}\}$ unless we are willing to assume some specific dependence structure for the sequence $\{(X_{1,j}, X_{2,j}) : j \geq 1\}$. If we assume, above and beyond (5) and (6), that

$$\{(X_{1,j}, X_{2,j}) : j \geq 1\} \text{ is an independent sequence of random vectors,} \tag{20}$$

then Theorem 1 of Majumdar and Majumdar (2017) shows (10) implies (13) (which renders Proposition 4 moot), by showing that the convergence of the Cesaro means of the sequence of cross sample correlation coefficients to 0 is a sufficient condition for (13) that turns out to be necessary for (19), equivalently, by Lemma 1, (10).

The assumption in (20), being the assumption in (8) with the identically distributed requirement removed, is weaker than (8). The aforesaid convergence of Cesaro means assumption is substantially weaker than (9). It is easy to see that if we combine the convergence of Cesaro means assumption with (5), (6), and (20), the resulting collection is weaker than the pair of assumptions in (8) and (9). All we have to do is to consider a pair of dependent but uncorrelated random variables and a sequence of iid copies of the resulting random vector. By Theorem 1 of Majumdar and Majumdar (2017), (13) can hold without (7). //

We now present the proofs of the results stated above.

Proof of Proposition 1 The key to the proof is the algebraic representation

$$W_{n_1,n_2} = \frac{\sqrt{n_1}(\bar{X}_1 - \mu_1)}{\sigma_1} \times V^{(1)}_{n_1,n_2} - \frac{\sqrt{n_2}(\bar{X}_2 - \mu_2)}{\sigma_2} \times V^{(2)}_{n_1,n_2}, \qquad (21)$$

where

$$\begin{aligned} V^{(1)}_{n_1,n_2} &= \frac{\sigma_1}{\sqrt{S_1^2 + \frac{n_1}{n_2} \times S_2^2}} \\ V^{(2)}_{n_1,n_2} &= \frac{\sigma_2}{\sqrt{\frac{n_2}{n_1} \times S_1^2 + S_2^2}}. \end{aligned} \qquad (22)$$

Now, let us fix $n_2$ and let $n_1 \to \infty$. Since $S_1^2$ converges (in probability) to $\sigma_1^2$ [Lemma A.5], $V^{(1)}_{n_1,n_2}$ converges to 0 and $V^{(2)}_{n_1,n_2}$ converges to $\sigma_2/S_2$. By the CLT and Slutsky's theorem, the first term in RHS(21) converges in distribution (and, by Theorem 4.2.9 of Fabian and Hannan (1985), in probability) to 0. Since the second term in RHS(21) converges in distribution to $\sqrt{n_2}(\bar{X}_2 - \mu_2)/S_2$, another application of Slutsky's theorem leads to the conclusion that, for fixed $n_2$, as $n_1 \to \infty$,

$$W_{n_1,n_2} \text{ converges in distribution to } -\sqrt{n_2}(\bar{X}_2 - \mu_2)/S_2. \qquad (23)$$

Since $S_2^2$ converges in probability to $\sigma_2^2$ [Lemma A.5], by the CLT and Slutsky's theorem, $-\sqrt{n_2}(\bar{X}_2 - \mu_2)/S_2$ converges, as $n_2 \to \infty$, in distribution to $Z_2 \sim \Phi$.

The same argument, with $n_1$ and $n_2$ interchanged, shows that if we fix $n_1$ and let $n_2 \to \infty$, $W_{n_1,n_2}$ converges in distribution to $\sqrt{n_1}(\bar{X}_1 - \mu_1)/S_1$, which, as $n_1 \to \infty$, converges in distribution to $Z_1 \sim \Phi$. □

Proof of Proposition 2 Clearly, $\{X_{i,j} : j \geq 1\}$ is an iid collection of standard Normal random variables, showing that (5) and (6) hold, with $\mu_i = 0$ and $\sigma_i = 1$. The measure induced by $D_j = X_{1,j} - X_{2,j}$ on $\mathcal{B}(\Re)$ is $\mathcal{N}_{2(1-\rho)}$ for every $j \in \mathfrak{N}$, implying that the measure induced by $\sum_{j=1}^{k} D_j/\sqrt{k}$ is $\mathcal{N}_{2(1-\rho)}$ for every $k \in \mathfrak{N}$. Since $S_{i,k}^2$ converges in probability to 1 as $k \to \infty$ [Lemma A.5] and $W_{\xi(k)} = \left(\sum_{j=1}^{k} D_j/\sqrt{k}\right)\left(S_{1,k}^2 + S_{2,k}^2\right)^{-1/2}$, by Slutsky's theorem the subnet $\{F_{\xi(k)} : k \in \mathfrak{N}\}$ of $\{F_\alpha : \alpha \in \mathfrak{N} \times \mathfrak{N}\}$ converges to $\mathcal{N}_{(1-\rho)}$, implying, by Lemma A.2 and the assumption $\rho \neq 0$, that (10) does not hold. □

Proof of Proposition 3 By Lemma A.3, to show (10) it suffices to show that given an arbitrary subnet $\{F_{\phi(\beta)} : \beta \in \mathcal{F}\}$ of $\{F_\alpha : \alpha \in \mathfrak{N} \times \mathfrak{N}\}$, there exists a further subnet $\{F_{\phi(\varphi(\delta))} : \delta \in \mathfrak{D}\}$ such that $\lim_\delta F_{\phi(\varphi(\delta))} = \Phi$. For $\alpha \in \mathfrak{N} \times \mathfrak{N}$, (22) implies

$$V_\alpha^{(i)} = \frac{\sigma_i}{\sqrt{e_\alpha^{1-i} S_{1,n_1(\alpha)}^2 + e_\alpha^{2-i} S_{2,n_2(\alpha)}^2}}. \tag{24}$$

Let $[0, \infty]$ denote the one-point compactification of $[0, \infty)$ [Dudley (1989, Theorem 2.8.1)]. Since every net taking values in a compact set has a convergent subnet [Lemma A.4], every subnet $\{e_{\phi(\beta)} : \beta \in \mathcal{F}\}$ of $\{e_\alpha : \alpha \in \mathfrak{N} \times \mathfrak{N}\}$ has a further subnet $\{e_{\phi(\varphi(\delta))} : \delta \in \mathfrak{D}\}$ such that

$$\lim_\delta e_{\phi(\varphi(\delta))} = \kappa \in [0, \infty]. \tag{25}$$

Since convergence in probability on Euclidian spaces is metrizable [Dudley (1989, Theorem 9.2.2)] and $n_i$ is order preserving and cofinal, by Lemma A.5, $S_{i,n_i(\alpha)}^2$ converges in probability to $\sigma_i^2$; consequently, by Lemma A.2, (24), and (25), in probability

$$\lim_\delta V_{\phi(\varphi(\delta))}^{(i)} = a_i = \frac{\sigma_i}{\sqrt{\kappa^{1-i}\sigma_1^2 + \kappa^{2-i}\sigma_2^2}} \in [0, 1]. \tag{26}$$

Note that $a_i$ depends on the subnet $\{e_{\phi(\varphi(\delta))} : \delta \in \mathfrak{D}\}$ through $\kappa$; $\kappa = 0$ implies $a_1 = 1$ and $a_2 = 0$, $\kappa = \infty$ implies $a_1 = 0$ and $a_2 = 1$, and, in general,

$$a_1^2 + a_2^2 = 1. \tag{27}$$

By (12) and Lemma A.2, $\{H_{\phi(\varphi(\delta))} : \delta \in \mathfrak{D}\}$ converges to $Q$. Since $W_\alpha = \langle Y_\alpha, V_\alpha \rangle$ by (21), where $V_\alpha = \left(V_\alpha^{(1)}, -V_\alpha^{(2)}\right)$, by Slutsky's theorem for nets [Lemma A.6] and Lemma A.7, (26) implies

$$\lim_\delta F_{\phi(\varphi(\delta))} = Q \circ \langle \cdot, (a_1, -a_2) \rangle^{-1}. \tag{28}$$

Since given any two unit vectors there exists an isometry mapping one to the other, by (26), (27), and the spherical symmetry of $Q$, RHS(28) does not depend on the subnet $\{e_{\phi(\varphi(\delta))} : \delta \in \mathfrak{D}\}$. By Lemma A.3, $\lim_\alpha F_\alpha = F \in \mathcal{M}(\mathfrak{R})$ exists for some $F \in \mathcal{M}(\mathfrak{R})$. Since an iterated limit exists and equals the double limit if the latter and the inner limit of the former exists [Lemma A.1], by (23), $\lim_{n_2 \to \infty} \lim_{n_1 \to \infty} F_{n_1, n_2} = F$; since $F = \Phi$ by Proposition 1, (10) follows. □

**Remark 7** DasGupta (2008, page 403) considers the Behrens-Fisher problem of comparing the means of two independent heteroscedastic normal populations and the two sample t-statistic that uses the pooled variance

$$S^2_{p,n_1,n_2} = \frac{(n_1-1)S_1^2 + (n_2-1)S_2^2}{n_1+n_2-2}$$

for studentization, that is

$$T_{n_1,n_2} = \frac{\overline{X}_1 - \overline{X}_2 - (\mu_1 - \mu_2)}{\sqrt{S_p^2\left(\frac{1}{n_1} + \frac{1}{n_2}\right)}},$$

as a potential pivot. He observes that if the ratio $n_1/n_2$ converges to 1, that is, the design is asymptotically balanced, then the asymptotic distribution of $T_{n_1,n_2}$ is standard Normal. As outlined below, that observation is a consequence of (10).

Given an arbitrary subnet $\{T_{\phi(\beta)} : \beta \in \mathcal{F}\}$ of $\{T_\alpha : \alpha \in \mathfrak{N} \times \mathfrak{N}\}$, there exists a further subnet $\{T_{\phi(\varphi(\delta))} : \delta \in \mathfrak{D}\}$ such that (25) holds. Since $n_i$ is order preserving and cofinal, and convergence in probability on Euclidian spaces is metrizable, by Lemma A.5 and Lemma A.2, in probability,

$$\lim_\delta S^2_{p,\phi(\varphi(\delta))} = \sigma_1^2\left(\frac{\kappa}{\kappa+1}\right) + \sigma_2^2\left(\frac{1}{\kappa+1}\right) = \frac{\kappa\sigma_1^2 + \sigma_2^2}{\kappa+1};$$

with

$$R_\alpha = \frac{S^2_{1,n_1(\alpha)} + e_\alpha S^2_{2,n_2(\alpha)}}{S^2_{p,\alpha}(1+e_\alpha)},$$

we obtain from Lemma A.5 and Lemma A.2 again, that, in probability,

$$\lim_\delta R_{\phi(\varphi(\delta))} = \frac{\sigma_1^2 + \kappa\sigma_2^2}{\kappa\sigma_1^2 + \sigma_2^2} = \theta(\kappa).$$

Since $T_\alpha = W_\alpha \sqrt{R_\alpha}$, with $K_\alpha \in \mathcal{M}(\mathfrak{R})$ denoting the distribution induced by $T_\alpha$, by Slutsky's theorem for nets, as long as (10) holds, we obtain

$$\lim_\delta K_{\phi(\varphi(\delta))} = \mathcal{N}_{\theta(\kappa)}.$$

For asymptotically balanced designs, $\kappa=1$ for every subnet $\{e_{\phi(\varphi(\delta))} : \delta \in \mathfrak{D}\}$, and Dasgupta's observation follows, showing that neither the Normality of the populations nor their cross sample independence is necessary for it. However, $T_{n_1,n_2}$ is an inferior choice for pivot compared to $W_{n_1,n_2}$, as asymptotic Normality of $W_{n_1,n_2}$ in (10) holds for all designs, whereas that of $T_{n_1,n_2}$ holds only for asymptotically balanced designs.   //

<u>Proof of Corollary 1</u> The proof of Proposition 3 applies verbatim once we replace $V_\alpha$ in that proof by $V_\alpha^* = \left(V_\alpha^{(1)}, V_\alpha^{(2)}\right)$ and note that $U_\alpha = \langle Y_\alpha, V_\alpha^* \rangle$. □

<u>Proof of Proposition 4</u> As observed in Remark 1, $L = \Phi \times \Phi$ if and only if $L \circ \langle \cdot, (a_1, a_2)\rangle^{-1} = \Phi$ for every unit vector $(a_1, a_2) \in \Re^2$.

Corresponding to every $a_1 \in [0, 1]$, there exist two unit vectors: $\left(a_1, \sqrt{1 - a_1^2}\right)$ in the first quadrant and $\left(a_1, -\sqrt{1 - a_1^2}\right)$ in the fourth quadrant. Since $\Phi$ is invariant under the map $x \mapsto -x$, it suffices to show that for every $a_1 \in [0, 1]$,

$$L \circ \left\langle \cdot, \left(a_1, \sqrt{1 - a_1^2}\right)\right\rangle^{-1} = \Phi = L \circ \left\langle \cdot, \left(a_1, -\sqrt{1 - a_1^2}\right)\right\rangle^{-1}. \tag{29}$$

We first show that given $\kappa \in [0, \infty]$, there exists a directed set $\mathcal{F}_\kappa$ and an order preserving and cofinal $\phi_\kappa : \mathcal{F}_\kappa \mapsto \mathfrak{N} \times \mathfrak{N}$ such that the subnet $\{e_{\phi_\kappa(\beta)} : \beta \in \mathcal{F}_\kappa\}$ of $\{e_\alpha : \alpha \in \mathfrak{N} \times \mathfrak{N}\}$ converges to $\kappa$.

If $\kappa = \infty$, let $\mathcal{F}_\infty = \mathfrak{N} \times \mathfrak{N}$ and define $\phi_\infty(\alpha) = \left(n_1(\alpha), \jmath\left(\sqrt{n_1(\alpha)}\right)\right)$, where $\jmath$ is the integer part function on $[1, \infty)$. Clearly, $\phi_\infty$ is order preserving. To show that $\phi_\infty$ is cofinal, given $\alpha \in \mathfrak{N} \times \mathfrak{N}$ choose $\widehat{\alpha} \in \mathfrak{N} \times \mathfrak{N}$ such that $n_1(\widehat{\alpha}) = \max\left(n_1(\alpha), (n_2(\alpha))^2\right)$; since $\jmath\left(\sqrt{n_1(\widehat{\alpha})}\right) \geq n_2(\alpha)$, $\phi_\infty(\widehat{\alpha}) \succeq \alpha$. Since $e_{\phi_\infty(\alpha)} = n_1(\alpha)/\jmath\left(\sqrt{n_1(\alpha)}\right) \geq \sqrt{n_1(\alpha)}$, the convergence of $e_{\phi_\infty(\alpha)}$ to $\kappa = \infty$ follows.

If $\kappa = 0$, let $\mathcal{F}_0 = \mathfrak{N} \times \mathfrak{N}$ and define $\phi_0(\alpha) = \left(n_2(\alpha), (n_2(\alpha))^2\right)$. Clearly, $\phi_0$ is order preserving. To see that $\phi_0$ is cofinal, given $\alpha \in \mathfrak{N} \times \mathfrak{N}$ choose $\widehat{\alpha} \in \mathfrak{N} \times \mathfrak{N}$ such that $n_2(\widehat{\alpha}) = \max(n_1(\alpha), n_2(\alpha))$; trivially, $(n_2(\widehat{\alpha}))^2 \geq n_2(\alpha)$, implying $\phi_0(\widehat{\alpha}) \succeq \alpha$. Since $e_{\phi_0(\alpha)} = n_2(\alpha)/(n_2(\alpha))^2 = 1/n_2(\alpha)$, the convergence of $e_{\phi_0(\alpha)}$ to $\kappa = 0$ follows.

If $\kappa \in (0, \infty)$, let $\mathcal{F}_\kappa = \mathfrak{N}$; define $\phi'_\kappa(r) = \jmath(\max((r+1)(\kappa - 1/r), 0) + 1)$ and $\phi_\kappa(r) = (\phi'_\kappa(r), r+1)$. It is easy to see that $\phi'_\kappa(r)$ is nondecreasing, implying that $\phi_\kappa$ is order preserving. To see that $\phi_\kappa$ is cofinal, given $\alpha \in \mathfrak{N} \times \mathfrak{N}$ choose $\widehat{r} \in \mathfrak{N}$ to equal $\max\left(\jmath\left(\left|\overset{*}{r}(\kappa, \alpha)\right| + 1\right), n_2(\alpha)\right)$, where

$$\overset{*}{r}(\kappa, \alpha) = \frac{\sqrt{(\kappa - 1 - n_1(\alpha))^2 + 4\kappa} - (\kappa - 1 - n_1(\alpha))}{2\kappa};$$

since $\kappa > 0$, $(\kappa - 1 - n_1(\alpha))^2 + 4\kappa$ is positive and $\overset{*}{r}(\kappa, \alpha)$ is well-defined. Since $\overset{*}{r}(\kappa, \alpha)$ is the bigger root of the quadratic (in $r$) $\kappa r^2 + (\kappa - 1 - n_1(\alpha))r - 1$, every $r \in \mathfrak{N}$ greater than $\overset{*}{r}(\kappa, \alpha)$ satisfies the inequality $\kappa r^2 + (\kappa - 1 - n_1(\alpha))r - 1 > 0$, equivalently, the inequality $(r+1)(\kappa - 1/r) > n_1(\alpha)$. Since $\widehat{r} > \left|\overset{*}{r}(\kappa, \alpha)\right| \geq \overset{*}{r}(\kappa, \alpha)$, $\phi'_\kappa(\widehat{r}) > n_1(\alpha)$, implying $\phi_\kappa(\widehat{r}) \succeq \alpha$. Finally, from the definition of $\jmath$,

$$\phi'_\kappa(r) - 1 \leq \max((r+1)(\kappa - 1/r), 0) < \phi'_\kappa(r),$$

implying

$$\max((\kappa - 1/r), 0) < \frac{\phi'_\kappa(r)}{r+1} \leq \max((\kappa - 1/r), 0) + \frac{1}{r+1};$$

since $\kappa > 0$ and $\phi_\kappa(r) = (\phi'_\kappa(r), r+1)$, $\lim_r e_{\phi_\kappa(r)} = \kappa$.

Given $a_1 \in [0,1]$, let $\kappa(a_1) = (\sigma_1 \sigma_2^{-1})^2(a_1^{-2} - 1) \in [0, \infty]$, where $\kappa(0)$ is interpreted to be $\infty$. Let $\left\{e_{\phi_{\kappa(a_1)}(\beta)} : \beta \in \mathcal{F}\right\}$ be a subnet of $\{e_\alpha : \alpha \in \mathfrak{N} \times \mathfrak{N}\}$ that converges to $\kappa(a_1)$. That, as in (26), implies, in probability

$$\lim_\beta V^{(1)}_{\phi_{\kappa(a_1)}(\beta)} = \frac{\sigma_1}{\sqrt{\sigma_1^2 + \kappa(a_1)\sigma_2^2}} = a_1$$

$$\lim_\beta V^{(2)}_{\phi_{\kappa(a_1)}(\beta)} = \frac{\sigma_2}{\sqrt{\frac{\sigma_1^2}{\kappa(a_1)} + \sigma_2^2}} = \sqrt{1 - a_1^2}.$$

Using $W_\alpha = \langle Y_\alpha, V_\alpha \rangle$ and $U_\alpha = \langle Y_\alpha, V_\alpha^* \rangle$ from the proofs of Proposition 3 and Corollary 1, (17), Lemma A.6, and Lemma A.7, (29) follows from (16) and (10). □

<u>Proof of Lemma 1</u> For $\alpha \in \mathfrak{N} \times \mathfrak{N}$, let

$$Q_\alpha = \frac{\frac{\sigma_1^2}{n_1(\alpha)} + \frac{\sigma_2^2}{n_2(\alpha)}}{\frac{S^2_{1,n_1(\alpha)}}{n_1(\alpha)} + \frac{S^2_{2,n_2(\alpha)}}{n_2(\alpha)}} = \frac{\sigma_1^2 + e_\alpha \sigma_2^2}{S^2_{1,n_1(\alpha)} + e_\alpha S^2_{2,n_2(\alpha)}} = \frac{e_\alpha^{-1}\sigma_1^2 + \sigma_2^2}{e_\alpha^{-1}S^2_{1,n_1(\alpha)} + S^2_{2,n_2(\alpha)}}. \quad (30)$$

By Lemmas A.6 and A.7, it suffices to show that $Q_\alpha$ converges in distribution to 1. Since convergence in distribution is metrizable, by Lemma A.3 it suffices to show that given an arbitrary subnet $\{Q_{\phi(\beta)} : \beta \in \mathcal{F}\}$ of $\{Q_\alpha : \alpha \in \mathfrak{N} \times \mathfrak{N}\}$, there exists a further subnet $\{Q_{\phi(\varphi(\delta))} : \delta \in \mathfrak{D}\}$ such that $Q_{\phi(\varphi(\delta))}$ converges in distribution to 1. Recall from (25) the existence of a subnet $\{e_{\phi(\varphi(\delta))} : \delta \in \mathfrak{D}\}$ that converges to $\kappa$. Since $\{e_{\phi(\varphi(\delta))} : \delta \in \mathfrak{D}\}$ converges in distribution to $\kappa$ and $S^2_{i,n_i(\alpha)}$ converges in probability (and hence, in distribution) to $\sigma_i^2$, the assertion follows from (30) (for $\kappa = \infty$, use the last representation) and Lemmas A.6 and A.7. □

*A.1. Double and iterated limits.* Given a metric space $(S, d)$, a $S$–valued double sequence is defined to be a function $x : \mathfrak{N} \times \mathfrak{N} \mapsto S$, where $\mathfrak{N}$ denotes the set of positive integers. We write $x(p,q)$ as $x_{p,q}$; recall that $x_{p,q}$ converges to $x \in S$ as $p, q \to \infty$ if, for every $\epsilon > 0$, there exists $n_0(\epsilon) \in \mathfrak{N}$ such that $p > n_0$ and $q > n_0$ imply $d(x_{p,q}, x) < \epsilon$.

Now suppose that for every fixed value of $p$, $\lim_{q \to \infty} x_{p,q} = y_p$ exists. Then, $\{y_p : p \in \mathfrak{N}\}$ is a $S$–valued sequence. If $\lim_{p \to \infty} y_p = y$ exists, then $y$ is an *iterated limit* of the double sequence $\{x_{p,q}\}$ and we write $\lim_{p \to \infty} \lim_{q \to \infty} x_{p,q} = y$. As illustrated in Section 8.20 of Apostol (1974), the existence of one iterated limit does not imply the existence of the other one (with $S = \mathfrak{R}$, let $x(p,q) = (-1)^p p(p+q)^{-1}$), the existence of both iterated limits does not imply their equality (with $S = \Re$, let $x(p,q) = p(p+q)^{-1}$), and the equality of the two iterated limits does not imply the existence of the double limit (with $S = \Re$, let $x(p,q) = pq(p^2+q^2)^{-1}$). However, an iterated limit exists and equals the double limit if the latter and the inner limit of the former exists.

**Lemma A.1** If the double limit of $\{x_{p,q}\}$, as $p,q \to \infty$, exists and is equal to $x$, then existence of $\lim_{q \to \infty} x_{p,q}$ for each fixed $p$ implies $\lim_{p \to \infty} \lim_{q \to \infty} x_{p,q} = x$.

Proof of Lemma A.1 This is Theorem 8.39 of Apostol (1974), stated and proved when $S$ is the complex plane with $d(a,b)$ being the absolute value of the difference $a - b$. That proof can be repeated verbatim here. □

*A.2. Nets and subnets.* A set $\mathfrak{D}$ endowed with a reflexive, anti symmetric, and transitive binary relation $\succeq$ is called a partially ordered set. The pair $(\mathfrak{D}, \succeq)$ is called a *directed set* if, for each $\beta, \gamma \in \mathfrak{D}$, there exists $\eta \in \mathfrak{D}$ such that $\eta \succeq \beta$ and $\eta \succeq \gamma$.

Given a metric space $(S, d)$ and a directed set $(\mathfrak{D}, \succeq)$, a $S$–valued net is defined to be a function $x : \mathfrak{D} \mapsto S$; we write the net as $\{x_\beta : \beta \in \mathfrak{D}\}$. Recall that the net $\{x_\beta : \beta \in \mathfrak{D}\}$ converges to $x \in S$ if, for every $\epsilon > 0$, there exists $\beta_0(\epsilon) \in \mathfrak{D}$ such that $\beta \succeq \beta_0(\epsilon)$ implies $d(x_\beta, x) < \epsilon$. It is worth recalling here that a $S$–valued sequence is a particular $S$–valued net.

Let $(\mathfrak{D}, \succeq)$ and $(\mathfrak{E}, \gg)$ be directed sets. Let $\phi : \mathfrak{E} \mapsto \mathfrak{D}$ be *order preserving*, that is, $i \gg j \Rightarrow \phi(i) \succeq \phi(j)$, and *cofinal*, that is, for each $\beta \in \mathfrak{D}$, there exists $\gamma \in \mathfrak{E}$ such that $\phi(\gamma) \succeq \beta$. Then the composite function $y = x \circ \phi$, where $x : \mathfrak{D} \mapsto S$, defines a net $\{y_\gamma : \gamma \in \mathfrak{E}\}$ in $S$, is called a subnet of $\{x_\beta : \beta \in \mathfrak{D}\}$, and is written as $\{x_{\phi(\gamma)} : \gamma \in \mathfrak{E}\}$.

**Lemma A.2** Let $\mathfrak{D}$ be a directed set and $\{x_\beta : \beta \in \mathfrak{D}\}$ a net taking values in $S$ that converges to $x \in S$. Then every subnet of $\{x_\beta : \beta \in \mathfrak{D}\}$ converges to $x$.

**Lemma A.3** Let $(\mathfrak{D}, \succeq)$ be a directed set and $\{x_\beta : \beta \in \mathfrak{D}\}$ a net taking values in $S$. Then $\{x_\beta : \beta \in \mathfrak{D}\}$ converges to $x \in S$ if and only if every subnet of $\{x_\beta : \beta \in \mathfrak{D}\}$ has a further subnet that converges to $x$.

**Lemma A.4** $S$ is compact if and only if every net in $S$ has a convergent subnet.

We refer the reader to Appendix A.1 of Majumdar and Majumdar (2017) for the proofs of these three lemmas.

*A.3. Miscellaneous results from probability.* We have used Lemmas A.5, A.6. A.7, and A.8 of this subsection in the paper.

**Lemma A.5 (Consistency of sample standard deviation)** As $k \to \infty$, $S_{i,k}$ converges almost surely (and hence, in probability) to $\sigma_i$.

Proof of Lemma A.5 Since

$$(k-1)S_{i,k}^2 = \sum_{j=1}^{k} X_{ij}^2 - k\bar{X}_{i,k}^2,$$

almost sure convergence follows from the SLLN. Since almost sure convergence of a sequence (as opposed to a net) of random variables implies convergence in probability [Fabian and Hannan (1985, Theorem 4.2.7)], the proof follows. □

**Lemma A.6 (Slutsky's theorem for nets)** Let $\mathcal{D}$ and $\mathcal{E}$ be metric spaces and $\mathfrak{T}$ a directed set. Let $\{X_\gamma : \gamma \in \mathfrak{T}\}$ and $\{Y_\gamma : \gamma \in \mathfrak{T}\}$ be nets of random elements taking values in $\mathcal{D}$ and $\mathcal{E}$, respectively, such that $X_\gamma \to X$ weakly and $Y_\gamma \to c$ weakly, where $X$ is a separable random element, that is, there exists a Borel measurable separable subset $\mathcal{M}$ of $\mathcal{D}$ such that $P(X \in \mathcal{M}) = 1$, and $c$ is a constant. Then, as a net of random elements taking values in $\mathcal{D} \times \mathcal{E}$, $(X_\gamma, Y_\gamma) \to (X, c)$ weakly.

Proof of Lemma A.6 See van der Vaart and Wellner (1996, Section 1.4). □

**Lemma A.7 (Continuous mapping theorem for nets)** Let $\{X_\gamma : \gamma \in \mathfrak{T}\}$ be a net of random elements taking values in $\Re^k$ such that $X_\gamma \to X$ in distribution. Let $g : \Re^k \mapsto \Re^m$ be a continuous function. Then $g(X_\gamma) \to g(X)$ in distribution in $\Re^m$.

Proof of Lemma A.7 See Theorem 1.3.6 of van der Vaart and Wellner (1996). □

**Lemma A.8** Let $S_1$ and $S_2$ be separable metric spaces such that $\{\mu_\beta : \beta \in \mathcal{F}\} \subset \mathcal{M}(S_1)$ and $\{\nu_\beta : \beta \in \mathcal{F}\} \subset \mathcal{M}(S_2)$ are two nets converging weakly to $\mu \in \mathcal{M}(S_1)$ and $\nu \in \mathcal{M}(S_2)$, respectively. Then the net $\{\mu_\beta \times \nu_\beta : \beta \in \mathcal{F}\} \subset \mathcal{M}(S_1 \times S_2)$ of product measures converges weakly to $\mu \times \nu$.

Proof of Lemma A.8 This lemma extends Lemma III.1.1 of Parthasarathy (1967) from sequences to nets. The proof of that lemma consists of three steps. In the first step it is shown that the algebra generated by all real valued functions on the product space that are products of bounded real valued uniformly continuous functions on the component spaces is dense (in the supremum norm) in the space of bounded real valued uniformly continuous functions on the product space. This step does not need any extension. In the second step it is shown (using the structure of this dense algebra) that the sequence of integrals of every function in that algebra wrt the sequence of product measures converges to the integral of the same function wrt the product of the limiting measures; this step readily extends to nets of measures from the definition of weak convergence of the component measures, thereby establishing (via the supremum norm approximation) that the net of integrals of every bounded real valued uniformly continuous function wrt the net of product measures converges to the integral of the same function wrt the product of the limiting measures. The third step uses Theorem II.6.1 of Parthasarathy (1967), which asserts the equivalence of weak convergence of a net of measures to another measure and the convergence of integrals of every bounded real valued uniformly continuous function wrt the net of measures to the integral of the same function wrt the limiting measure, so no extension is required for this step. □

RAJESHWARI MAJUMDAR
rajeshwari.majumdar@uconn.edu
PO Box 47
Coventry, CT 06238

SUMAN MAJUMDAR
suman.majumdar@uconn.edu
1 University Place
Stamford, CT 06901